\documentclass[12pt]{amsart}
\usepackage{amssymb,amsmath,amsfonts,amsthm,amscd,enumerate}

\usepackage{graphics}
\usepackage[colorlinks=true]{hyperref}
\usepackage{amssymb,amsmath,amsfonts,amsthm,amscd,mathrsfs}

\input xy
\xyoption{all}
\newcommand{\CC}{\mathbb C}
\newcommand{\RR}{\mathbb R}

\newcommand{\maC}{\mathcal{C}}
\newcommand{\maD}{\mathcal{D}}

\newcommand{\ZZ}{\mathbb{Z}}
\newcommand{\Sch}{\mathscr{S}}
\newcommand{\NN}{\mathbb{N}}

\newcommand{\maE}{\mathcal{E}}

\newcommand{\maG}{\mathcal{G}}
\newcommand{\maN}{\mathcal{N}}

\newcommand{\ep}{\varepsilon}

\newcommand{\smooth}{\mathcal{C}^{\infty}}

\newcommand{\ip}[1]{\langle #1 \rangle}

\newcommand{\gp}[1]{\Gamma^{\infty}(#1)}
\theoremstyle{definition}
\newtheorem{definition}{Definition}[section]

\theoremstyle{plain}
\newtheorem{theo}[definition]{Theorem}
\newtheorem{prop}[definition]{Proposition}
\newtheorem{lem}[definition]{Lemma}
\newtheorem{cor}[definition]{Corollary}

\theoremstyle{remark}
\newtheorem{remark}[definition]{Remark}
\newtheorem{eg}[definition]{Example}

\begin{document}

\textcolor{red}{
 \title[Geometrical embeddings of distributions into $\maG$]{Geometrical embeddings of 
distributions into algebras of generalized functions} 
}
\author{Shantanu Dave}
\address{University of Vienna, Faculty of Mathematics, Nordbergstrasse 15, 1090 Vienna, Austria}
\email{shantanu.dave@univie.ac.at}
\thanks{Supported by FWF grant Y237-N13 of the Austrian Science Fund.}

\maketitle
\begin{abstract}
We use spectral theory to produce embeddings of distributions into algebras of generalized functions 
on a closed Riemannian manifold. These embeddings are invariant under isometries and  preserve 
the singularity structure of the distributions.
\end{abstract}
%\tableofcontents
\section{Introduction}
Algebras of generalized functions as introduced by Colombeau \cite{c1,c2,Biag,MObook,book} provide a 
consistent framework for carrying out nonlinear analysis on spaces of distributions, 
especially in nonlinear partial differential equations and regularity theory 
\cite{MObook,garetto_hoermann,guezygmund,HOP,delcroix}. There are various choices of embedding distributions 
into these algebras, representing fundamental aspects of nonlinear modeling depending 
on the choice of regularizing process.

The special algebra of generalized functions \cite{c1,book} over an open set in $\RR^n$ or on a manifold, 
is characterized in terms of a regularizing net of smooth functions and asymptotic estimates in terms 
of the regularizing parameter $\varepsilon$. More precisely, let $I=(0,1]$ then the special algebra 
over $\RR^n$ is defined by the following quotient construction. Firstly the moderate nets are defined as
\begin{align*}
\maE_M(\RR^n)=\{&(U_{\varepsilon})_{\varepsilon\in I}\in\smooth(\RR^n)^I|\forall\, 
\textrm{semi-norms}\,\rho\,\textrm{on}\,\smooth(\RR^n)\\
&\exists N\in\ZZ \;\textrm{such that}\; \rho(U_{\varepsilon})\sim O(\varepsilon^N)\}.
\end{align*}
Here the notation $f(\varepsilon)\sim O(g(\varepsilon))$ implies that there exists 
$\varepsilon_0>0$ and a constant $C>0$ such that   $|f(\varepsilon)|<Cg(\varepsilon)$ for all 
$\varepsilon<\varepsilon_0\,$.

The nets of negligible growth are given by
\begin{align*}
\maN(\RR^n)=\{&(U_{\varepsilon})_{\varepsilon\in I}\in\smooth(\RR^n)^I|\forall\, 
\textrm{semi-norms}\,\rho\,\textrm{on}\,\smooth(\RR^n)\\
&\forall  N\in\ZZ \,: \rho(U_{\varepsilon})\sim O(\varepsilon^N)\}.
\end{align*}
form an ideal in $\maE_M(\RR^n)$ and the special Colombeau algebra is defined as the quotient
\begin{align*}
&\maG(\RR^n):=\maE_M(\RR^n)/\maN(\RR^n).
\end{align*}
 Given a moderate net of smooth functions $(u_{\varepsilon})_{\varepsilon}$ one  usually denotes the corresponding element in the quotient  by $[(u_{\varepsilon})_{\varepsilon}]$ or  by $\ip{u_{\varepsilon}}$. We shall often identify the net with the element it represents if no confusion can arise. Also, for the sake of simplicity we shall omit the index
$M$ in the notation of moderate nets $\maE_M$ henceforth and simply write $\maE$. The algebra $\smooth(\RR^n)$ of smooth functions can of course be embedded into $\maG(\RR^n)$ as a subalgebra by the constant embedding $\sigma:\,\maC^\infty(\RR^n)\rightarrow \mathcal G(\RR^n)$, $f\mapsto [(f)_{\varepsilon}]$. To embed nonsmooth  distributions one first picks a mollifier $\rho\in\Sch(\RR^n)$ such that the  net $\rho_{\varepsilon}:=\varepsilon^{-n}\rho(x/\varepsilon)$ suitably approximates the delta distribution. The compactly supported distributions are embedded by $\omega\rightarrow [(\omega*\rho_{\varepsilon})_{\varepsilon}]$. This embedding can be extended to all distributions in $\maD'(\RR^n)$ by suitable partitions of unity and cut off functions using the sheaf structure. The resulting embedding $\iota_{\rho}:\maD'(\RR^n)\rightarrow \maG(\RR^n)$ commutes with partial derivations and matches up with the constant embedding $\sigma$ of smooth functions mentioned above. This implies that the inclusion of smooth functions via $\iota_{\rho}$ is in fact an algebra homomorphism. 

Taking a closer look at the above process one realizes that:
\begin{enumerate}[a{)}]
\item Convolution by the smooth mollifier is in fact a smoothing operator, that is it maps tempered distributions to smooth function.
\item Since the net $\rho_{\varepsilon}\rightarrow \delta$, the convolution by this net gives a net of smoothing operators that approximate the identity operator (in a sense to be made precise below).
\end{enumerate}
  In view of the above observations we consider embedding distributions into algebras of generalized functions using a net of smoothing operators approximating the identity operator. In case of a compact manifold one could rely on operators defined either by smooth kernels or by functional calculus to generate such nets of smoothing operators.
  In this paper we shall work with a compact Riemannian manifold $M$ with associated Laplace operator $\Delta$.  We recall that by functional calculus any Schwartz class function $F\in\Sch(\RR)$ defines a smoothing operator $F(\Delta)$. Let us set $F_{\varepsilon}(x)=F(\varepsilon x)$. The main result of this paper is then as follows:
\begin{theo}\label{main}
If a Schwartz function $F$ is $1$ near the origin then the net of smoothing operators $F_{\varepsilon}(\Delta)$ provides an embedding $F_{\varepsilon}:\maD'(M)\rightarrow \maG(M)$. Such an embedding preserves the multiplication on $\smooth(M)$ (coincides with the constant embedding for smooth functions) and is invariant under isometries.
\end{theo} 
 We emphasize the following interesting features of our approach.
 The embeddings $F_{\varepsilon}$ are of course globally constructed. Thus we do not need to work with local coordinates and use the sheaf structure of the algebra $\maG$. In addition the global construction makes it 
compatible with the geometry of the Riemannian manifold as these embeddings  commute with the Laplace operator and hence are preserved under isometries.
The embeddings  introduced so far in the literature (\cite{AB,RD,ndg})
are `non-geometric' in that they depend on choices of partitions of unity, cut-offs,
etc. (cf.\ \cite{book}, Sec.\ 3.2.2 for a discussion).

 We also note that this functional calculus technique can be naturally applied in case of embedding distributional sections of a (complex) vector bundle into the 
corresponding space of generalized sections.
%$\tilde{\CC}$ module of locally convex space defined by the smooth sections.
 
The proof of Theorem \ref{main} is based on the functoriality of the construction of generalized 
$\tilde{\CC}$-modules over a locally convex space. The negligibility and moderateness estimate are hence obtained as direct consequences of Weyl's estimates for eigenvalues of positive elliptic operators.  
For instance the Weyl estimates being asymptotic provides the well known continuity of the map 
\begin{eqnarray*}
\Sch(\RR) &\rightarrow& \Psi^{\!-\!\infty}(M)\\
f &\rightarrow& f(\Delta)  
\end{eqnarray*}
and hence establish the moderateness estimates.
\section{Locally convex spaces} \label{functor}
Let $X$ be a locally convex (Hausdorff) topological vector space then one can associate
a generalized locally convex space $\maG_X$ (see \cite{Garetto,garetto_duals}, as well as
\cite{Scarpal}) as follows.
Let $I\subset \RR$ be the interval $(0,1]$. Define the moderate nets in $X^I$ to
be, 
\begin{align*}
\maE(X):= \{(x_{\varepsilon}): \, \forall\,&\textrm{continuous semi-norm}\,\rho \textrm{ on } X\,\, \exists
\textrm{ an integer}\,N\\
&\,\textrm{such that}\,\rho(x_{\varepsilon})\sim O(\varepsilon^N)  
\}.
\end{align*}
 Similarly we can define the negligible nets to be
\begin{align*}
\maN(X):= &\{(x_{\varepsilon}): \, \forall\,\textrm{continuous semi-norm}\,\rho \,\, 
 \textrm{ on} \,\,  X \textrm{ and}\, \forall m\\
&|\rho(x_{\varepsilon})|\sim O(\varepsilon^m) \}.
\end{align*}
The generalized locally convex space of $X$ is then defined to be the quotient,
\[
\maG_X:=\maE(X)/\maN(X).
\]
One notes that in defining $\maE(X),\,\maN(X)$ it suffices to restrict to a family of
semi-norms that generate the locally convex topology for $X$. 
When $X=\smooth(M)$ is the space of smooth functions on a manifold $M$ then we
also write $\maG(M):=\maG_{\smooth(M)}$. Also for $X=\CC$ we call $\maG_{\CC}$ the space of generalized numbers and denote it by $\tilde{\CC}$. $\tilde{\CC}$ is the ring of constants in $\maG_{\CC}$. Every $\maG_X$ is naturally a module over the ring $\tilde{\CC}$, and hence  is often referred to as the $\tilde{\CC}$-module associated with $X$.  The sharp topology on $\tilde{\CC}$ is the topology generated by sets of the form $U_{x,p}$ where $x\in\tilde{\CC}$, $p$ is an integer and
\[U_{x,p}:=\{\ip{x_{\varepsilon}-\varepsilon^p,x_{\varepsilon}+\varepsilon^p}| \ip{x_{\varepsilon}}=x\}.\]
Any continuous semi-norm $\rho$ on a locally convex spaces $X$  by definition provides a  map $\tilde{\rho}:\maE(X)\rightarrow \maE(\CC)$ by applying $\rho$ to each component. In fact $\tilde{\rho}$ descends to a map 
$\maG_X\rightarrow \tilde{\CC}$.
The sharp topology on any $\tilde{\CC}$-module $\maG_X$ shall be defined as the weakest topology that 
makes all the above $\tilde{\rho}$'s continuous.

We recall the functoriality of the above construction.
\begin{lem}\label{lcs}
If $\phi:X\rightarrow Y$ is a continuous linear map between locally
convex spaces $X $ and $Y$ then there is a natural induced map $\phi_*:\maG_X\rightarrow
\maG_Y$ defined on the representatives as
$\phi_*((x_{\varepsilon}))=(\phi(x_\varepsilon))$. 
Further $\phi_*$ is continuous with respect to the sharp topology.
\end{lem}
\begin{proof}
If $\tau$ is a continuous semi-norm on $Y$ then $\tau\circ \phi$ is a continuous semi-norm on $X$. Thus if $x_{\varepsilon}$ satisfies an estimate with $\tau\circ \phi$ in $\maE(X)$ or $\maN(X)$ then $\phi(x_{\varepsilon})$ satisfies the exact same estimates with respect to $\tau$ in $\maE(Y)$ or $\maN(Y)$. Thus $\phi_*$ is well-defined.
 Since basic open sets $U$ in $\maG_Y$ are pull-backs of open sets in $\tilde{\CC}$ by some semi-norm $\rho$, thus  $\phi^{-1}_*(U)$ is a pullback of an open set with respect to $\rho\circ \phi$.
\end{proof}
For example any smooth map between two manifold $f:M\rightarrow N$ gives rise
to a pull back map $f^*:\maG(N)\rightarrow \maG(M)$.
\section{Embedding of distributions}
 Let $M$ be a compact manifold without boundary. A continuous linear operator $T:\maD'(M)\rightarrow \smooth(M)$ is
 called a smoothing operator.  We shall denote the space of all smoothing
 operators by $\Psi^{-\infty}(M)$. Each smoothing operator is a pseudodifferential operator and extends a map $\hat{T}:\smooth(M)\rightarrow \smooth(M)$.
 The following is a well know characterization of smoothing  operators:
 %\begin{lem}
 Let $dy$ be a Riemannian  density on $M$. An operator $T$ is smoothing if and only if there
 is a smooth function $k(x,y)\in \smooth(M\times M)$ such that for any
 $u\in\smooth(M)$,
 \begin{equation}\label{ker}
Tu(x)=\int_M k(x,y)u(y)dy.
\end{equation}
It also follows that any  such identification gives an isomorphism of $\Psi^{-\infty}(M)$ with  $\smooth(M\times M)$ as locally convex spaces.
 In the sequel we shall always assume a Riemannian  density  is provided to us. This in particular shall imply that we are given an embedding of $\smooth(M)$ into $\maD'(M)$.
 \begin{definition}\label{abstract_embedding}
  A parametrized family $T_{\varepsilon}$ of smoothing operator is
 called a special embedding if
 \begin{enumerate} 
 \item \label{moderate}For any semi-norm $\rho$ on $\smooth(M)$ and any distribution $u$ there
 exists an integer $N$ such that,
\[
\rho(T_{\varepsilon}u)\sim O(\varepsilon^N).
\]
\item \label{embed}
Let $u\in \maD'(M)$ be such that $\rho(T_{\varepsilon}u)\sim O(\varepsilon^N)$ for all integers $N$ and all semi-norms $\rho$  on $\smooth(M)$, then $u=0$.

 \item\label{negligible} If $f\in\smooth(M)$ is a smooth function  on M then for all semi-norms $\rho$ on $\smooth(M)$ and given any integer $N$,
\[
\rho(T_{\varepsilon}f-f)\sim O(\varepsilon^N).
\]
 \end{enumerate}
\end{definition}
The above definition is tailored to obtain an embedding:
\begin{lem}
If a collection of smoothing operators $T_{\varepsilon}$ is a special embedding then it provides an embedding of
the space of distributions $\maD'(M)\rightarrow \maG(M)$ namely $u\mapsto [(T_{\varepsilon}u)_{\varepsilon}]$ which coincides with the constant embedding of smooth functions and hence preserves the
multiplication of smooth functions.
\end{lem}
%\begin{proof}
%By Definition \ref{abstract_embedding} \eqref{moderate}, the map $u\rightarrow T_{\varepsilon}u$ maps $\maD'(M) $ to $\maG(M)$.
%That the map $u\rightarrow \ip{T_{\varepsilon}u}$ to the quotient is injective follows from \ref{abstract_embedding} \eqref{embed}.
%And finally $\ip{T_{\varepsilon}f}=[(f)_{\varepsilon}]$ for a smooth function $f$ follows from \ref{abstract_embedding} \eqref{negligible}.
%\end{proof}
\begin{definition}
Let $k>0$ be an integer. We say that a parametrized family $T_{\varepsilon}$ is an
order-$k$ embedding if $T_{\varepsilon}$ satisfies \eqref{moderate} and \eqref{embed} above and
the condition that for any smooth function $f\in \smooth(M)$ we have
\[
\rho(T_{\varepsilon}f-f)\sim O(\varepsilon^k).
\]
 \end{definition}
\begin{eg}
Let $M=\RR^n$  and let $\rho\in\Sch(\RR^n)$ be in the Schwartz class such that
$\hat{\rho}(0)=1$ and all moments of $\rho$ vanish. Let $\rho_{\varepsilon}(x)=\frac{1}{\varepsilon ^n}\rho(x/\varepsilon)$ be the approximate
unit  converging to the delta function. The operators $T_{\varepsilon}(f):=\rho_{\varepsilon}\ast f$ form a one parameter family of operators
mapping the tempered distributions to smooth functions which satisfies all the
above properties in the class of tempered distributions and smooth tempered
functions. Of course it does not include all distributions on $\RR^n$ because of
non-compactness.
\end{eg}
Now  we fix a Riemannian metric on $M$ and let $\Delta$ be the corresponding
Laplace operator. Let $f\in \Sch(\RR)$ then by functional calculus the operator $f(\Delta)$ is a smoothing operator. We
shall use smoothing operators obtained from the geometric  Laplace operator to
obtain geometrical  embeddings.
\section{Geometric Laplace operators and Functional calculus}\label{funct_calc}
 Let $M$ be a closed Riemannian manifold. Let $\Delta$  be the associated scalar
 Laplace operator. Let $\phi_k$ in $L^2(M)$ be  eigenfunctions of $\Delta$ with
 eigenvalues $\lambda_k$. The classical theorem of Hermann Weyl implies that there exist $a,\,b>0$ such that
 $\frac{\lambda_n}{an^b}\rightarrow 1$ as $n\rightarrow\infty$. This together with  elliptic regularity  of
 $\Delta$ gives that the ``Fourier series'' expansion
\[f=\sum_ka_k\phi_k,\quad f\in\smooth(M),\] provides an
 isomorphism
 \[\overline{\Delta}: \smooth(M)\rightarrow \Sch(\NN)\quad \textrm{by}\quad f\mapsto
 (a_k).\]

In the following $\hat{\otimes}$ as usual denotes the projective tensor product.
\begin{prop}\label{projective}If $M$ and $N$ are closed manifolds then  
\[\smooth(M)\,\hat{\otimes}\,\smooth(N)\simeq \smooth(M\times N).\]
\end{prop}
\begin{proof} We give an outline of the proof here. See \cite{Grothendieck} for original and  more detailed discussion.  Fix a metric on both $M$ and $N$. Let $\Delta_M$ and $\Delta_N$ be Laplace operators on $M$ and $N$ respectively. Then the Laplace on $M\times N$  with the product metric is given by $\Delta_{M\times N}:=\Delta_M\otimes 1+1\otimes \Delta_N$. Then one can check that the following diagram commutes,
\begin{center}
$ \xymatrix{ {\smooth(M)\,\hat{\otimes}\,
\smooth(N)} \ar@<3ex>[d]^{\overline{\Delta}_M}  \ar@<-3ex>[d]_{\overline{\Delta}_N}\ar[r]&{\smooth(M\times N)}\ar[d]^{\overline{\Delta}_{M\times N}}\\
 {\Sch(\NN)\, \hat{\otimes}\, \Sch(\NN)}\ar[r]&{\Sch(\NN\times \NN)}}$
\end{center}
Note that the map on the top row is an isomorphism since the one on the lower row is.
\end{proof}
We view the smooth functions in $\smooth(M\times M)$ as kernels of smoothing
operators. Thus in particular the above proposition $\Psi^{-\infty}(M)$ gives the following .
\begin{cor}
Let $a_k$ be a sequence in $\Sch(\NN)$. Then the sum
$K(a_k)=\sum a_k\phi_k\otimes\phi_k$ is in $\smooth(M\times M)$ and hence defines a
smoothing operator. 
\end{cor}
For any Schwartz class function $f\in\Sch(\RR)$ we denote by $f(\Delta)$ the
operator $\sum_k f(\lambda_k)\phi_k\otimes \phi_k$. Since $f(\Delta) $ is
diagonal with the same eigenspaces as $\Delta$ it automatically commutes with
$\Delta$.

We thus obtain a map of
locally convex spaces $\Sch(\RR)\rightarrow \Psi^{-\infty}(M)$ namely $f\rightarrow f(\Delta)$  which, abusing notation, we also denote by $\Delta$.
\begin{lem}\label{func_calc}
The map $\Delta:\Sch(\RR)\rightarrow \Psi^{-\infty}(M)$ is continuous.
\end{lem}
\begin{proof}
Let $\lambda_n,\phi_n$ be a spectral decomposition for $\Delta$ as above. Then by the kernel theorem (\ref{ker}) and by Proposition \ref{projective} one can identify $\Psi^{_\infty}(M)$ with $\Sch(\NN\times \NN)$  isomorphically. The map $\Delta$ now becomes
\[\Sch(\RR)\ni f\rightarrow \{a_{ij}\}=\left\{\begin{array}{cl}
f(\lambda_i) \hskip 0.5in &\text{if }i=j\\
0 \hskip 0.5in
&\text{otherwise}
\end{array}.\right. \] 

The ideal $I_{sp(\Delta)}:=\{g\in\Sch(\RR)|g(\lambda_k)=0\;\textrm{for all}\; k\}$ is a closed ideal. Now the map $\Delta$  is a composition of two maps $q:\Sch(\RR)\rightarrow \Sch(\NN)$, the quotient map with respect to the ideal $I_{sp(\Delta)}$, and the diagonal embedding $ \Sch(\NN)\rightarrow \Sch(\NN\times \NN)$. Note that  the quotient $\Sch(\RR)/I_{sp(\Delta)}\simeq \Sch(\NN)$  since by Weyl's theorem the $\lambda_n$ grow moderately.
\end{proof}
As a consequence of Lemma \ref{lcs} we get an induced map
\begin{eqnarray}
\Delta_*:\maG_{\Sch(\RR)}\rightarrow \maG_{\Psi^{-\infty}(M)}.
\end{eqnarray}
The map $\Delta_*$  will provide us with estimates on $\maG_{\Psi^{-\infty}(M)}$ in terms of  
those in $\maG_{\Sch(\RR)}$.
 \section{Geometrical  Embeddings}
 Let $F\in\Sch(\RR)$ and $F\equiv 1$ near the origin in an interval $(-t,t)$.  
Set $F_{\varepsilon}(x):=F(\varepsilon x)$ for $0<\varepsilon<1$, then the $F_{\varepsilon}$
 form a net of approximate units in $\Sch(\RR)$. We fix one such net for the current section. Let $M$ be a closed Riemannian manifold and $\Delta$ the associated scalar
  Laplace operator. Further let $(\lambda_k,\phi_k)$ be eigenvalues and eigenvectors of
  $\Delta$, where $0\leq \lambda_1\leq\lambda_2\ldots$ are counted with
  multiplicity. We would like to analyze the  net $F_{\varepsilon}(\Delta)$ of smoothing
  operators. To begin with we set
\begin{eqnarray}\label{Weyl}  
N_{\varepsilon}:=\max\{n\,|\,\lambda_n<\frac{t}{\varepsilon}\}.
\end{eqnarray}
 \begin{lem}\label{asymptotic}
 There exist $C>0,\,\alpha>0$ and $\varepsilon_0$ such that $N_{\varepsilon}\geq C\varepsilon^{-\alpha}$
for all $\varepsilon<\varepsilon_0$.
\end{lem}
\begin{proof}
This is a direct consequence of Weyl's estimates which  provide an asymptotic estimate for the spectral counting function 
\[N_{\Delta}(\lambda)=\#\{\lambda_k|\,\lambda_k<\lambda\}.\] 

(See \cite{Berger, BGV,Rudin,Roe}.)

Let $m=\operatorname{dim}(M)$ then
\[N_{\Delta}(\lambda) \sim \frac{\operatorname{vol}(M)} {(4\pi)^{\frac{m}{2}}\Gamma(m/2+1)}\lambda^{\frac{m}{2}}.\]
Hence for $\alpha\geq\frac{m}{2}$ we have $N_{\varepsilon}\geq C\varepsilon^{-\alpha}$ for $\varepsilon$ small enough.
\end{proof}
\begin{lem}\label{base_estimate}
Let $f\in\smooth(M)$  then for any natural number  $N$,
\[ 
\|F_{\varepsilon}(\Delta)f-f\|_{L^2(M)}\sim O(\varepsilon^N).
\]
 \end{lem}
 \begin{proof}
 By Weyl's estimates as before the Fourier expansion in eigenfunctions
 \[f=\sum a_k\phi_k
\]
implies $(a_k)\in\Sch(\NN)$.
 Let $\alpha$ be as in Lemma \ref{asymptotic}. Given any integer $N$, we can find a $k$ such that for all $l>k$
 \[\sum_{n>l}|a_n|^2<\frac{1}{l^{\frac{N}{\alpha}}}.\]
 Now choose $\varepsilon_0>0$ such that $\lambda_k<\frac{t}{\varepsilon_0}$. Then for
 any $\varepsilon<\varepsilon_0$

\[F_{\varepsilon}(\Delta)f-f=\sum_{n>N_{\varepsilon}}a_n(F_{\varepsilon}(\lambda_n)-1)\phi_n.\]
 Since $F_{\varepsilon}(\lambda_n)$ are all uniformly bounded, 
 \begin{align*}
 \|F_{\varepsilon}(\Delta)f-f\|_{L^2(M)}^2&\leq C\sum_{n>N_{\varepsilon}}|a_n|^2\\
 &\leq \frac{1}{N_{\varepsilon}^{\frac{N}{\alpha}}}\sim O(\varepsilon^N).
\end{align*}
 \end{proof}
 Let $D:=(1+\Delta)^{\frac{1}{2}}$. Then for any integer $k$ the
 Sobolev space $H^k(M)$ is the completion of $\smooth(M)$ with respect to the
 norm 
 \[\|f\|_k:= \|\maD^kf\|_{L^2(M)}.\]

 Since each $F_{\varepsilon}(\Delta)$ is a smoothing operator it is a bounded
 operator $H^k(M)\rightarrow L^2(M)$ for any integer $k$. For a bounded operator 
 $T: H^k(M)\rightarrow L^2(M)$ denote by $\|T\|_k$ its operator norm. 
\begin{lem}\label{norm_estimate}
For any integer $k$ the map $T\rightarrow \|T\|_k$ is a continuous semi-norm on $\Psi^{-\infty}(M)$.
\end{lem} 
\begin{proof} Let  $k(x,y)$ be the smooth kernel of $T$. Let $u(y)$ be a smooth function, then
using Cauchy-Schwartz on the integral operator,
\begin{align*}
\|Tu\|_{L^2(M)}&=\left(\int_M|\int_M k(x,y)u(y)dy|^2dx\right)^{\frac{1}{2}}\\
&=\left(\int_M|\int_M (1+\Delta)^{-\frac{k}{2}}k(x,y)(1+\Delta)^{\frac{k}{2}}u(y)dy|^2dx\right)^{\frac{1}{2}}\\
&\leq \|u\|_k\operatorname{Vol}^{\frac{1}{2}}(M)\left(\int_M\int_M (1+\Delta)^{-\frac{k}{2}}|k(x,y)|dy^2dx\right)^{\frac{1}{2}}.
\end{align*}
Note that we have used the fact that $(1+\Delta)^{-s}$ is a self-adjoint operator on $L^2(M)$. Now $\smooth(M)$ is dense in all Sobolev spaces so the above norm estimates hold on $H^k(M)$ and can be obtained from the norm of the kernel in $\smooth(M\times M)$.
\end{proof}
\begin{lem}
 Let $F$ be a function in the Schwartz class $\Sch(\RR)$ and let $F_{\varepsilon}(x):=F(\varepsilon x)$. Then there exists an integer $q$ such that  $\|F_{\varepsilon}(\Delta)\|_k\sim O(\varepsilon^q)$.
 \end{lem}
 \begin{proof}
One notes that $F_{\varepsilon}(x)$ is a moderate net in $\Sch(\RR)$. Since $\Delta:\Sch(\RR)\rightarrow \Psi^{-\infty}(M)$ by Lemma \ref{func_calc} is a continuous map of locally convex vector spaces by Lemma \ref{lcs} $F_{\varepsilon}(\Delta)$ is a moderate net in $\Psi^{-\infty}(M)$. The result follows  by continuity of the semi-norm $\|\,\|_k$.
 \end{proof}
Now here are our first embeddings.
\begin{prop}\label{table_theorem}
Let $F$ be a Schwartz function with $F\equiv 1$ near $0$ then
$T_{\varepsilon}=F_{\varepsilon}(D)$ is a special embedding of $\maD'(M)$ into
$\maG_M$.
\end{prop}
 \begin{proof}
We must check the three conditions in the Definition \ref{abstract_embedding} of special embeddings.
 \begin{enumerate}
\item Since $\maD'(X)=\bigcup_s H^s(M)$ for any distribution $u$ there exists an $s$ such that $u\in H^s(M)$. Considering $F_{\varepsilon}(\Delta)$ as an operator from $H^s(M)\rightarrow L^2(M)$ one has:
\begin{align*}
\|F_{\varepsilon}(\Delta)u\|_{L^2(M)}&\leq \|F_{\varepsilon}(\Delta)\|_s||u\|_s\\
&\simeq O(\varepsilon^N)\qquad \textrm{by  Lemma \ref{norm_estimate}.}
\end{align*}
Now for any other Sobolev norm estimate one notices that $F_{\varepsilon}(\Delta)$  commutes with  $(1+\Delta)^{\frac{k}{2}}$ and applying Lemma \ref{norm_estimate} to  $(1+\Delta)^{\frac{k}{2}}u\in H^{k+s}(M)$ we obtain,
\begin{align*}
\|F_{\varepsilon}(\Delta)u\|_{H^k(M)}&= \|(1+\Delta)^{\frac{k}{2}}F_{\varepsilon}(\Delta)u\|_{L^2(M)}\\
&=\|F_{\varepsilon}(\Delta)(1+\Delta)^{\frac{k}{2}}u\|_{L^2(M)}\\
&=\|F_{\varepsilon}(\Delta)\|_{s+k}\|(1+\Delta)^{\frac{k}{2}}u\|_{s+k}\simeq O(\varepsilon^N).
\end{align*}
\item
For any $u\in \maD'(M)$ we shall show that $\lim_{\varepsilon\rightarrow 0}F_{\varepsilon}(\Delta)u =u$ therefore the map $\maD'(M)\rightarrow \maG(M)$ is injective. First when $u\in L^2(M)$ then we write a Fourier expansion $u=\sum a_n\phi_n$. And for $\varepsilon<\frac{t}{\lambda_k}$ one observes that $u-F_{\varepsilon}(\Delta)(u)\leq\sum_{n>k}a_n(F_{\varepsilon}(\lambda_n)-1)\phi_n$, where the right hand side tends to zero in norm as $k\rightarrow \infty$. For $u$ in any other Sobolev space $H^s(M)$ we notice that the above argument can be applied to $(1+\Delta)^{-\frac{s}{2}}u$ and that these operators commute with $F_{\varepsilon}(\Delta)$.
\item
By G\r{a}rding's inequality 
\[\|u\|_{H^{k+1}(M)}\leq C(\|u\|_{H^k(M)}+\|\Delta^{\frac{1}{2}}u\|_{H^k(M)},\]
hence by applying  Lemma \ref{base_estimate} one obtains $\forall N$,
\begin{align*}
\|F_{\varepsilon}(\Delta)f-f\|_1&\leq  \|F_{\varepsilon}(\Delta)f-f\|_{L^2(M)}+ \|\Delta^{\frac{1}{2}} (F_{\varepsilon}(\Delta)f-f)\|_{L^2(M)}\\
&\leq  \|F_{\varepsilon}(\Delta)f-f\|_{L^2(M)}+ \|F_{\varepsilon}(\Delta)\Delta^{\frac{1}{2}}f-\Delta^{\frac{1}{2}}f\|_{L^2(M)}\\
&\sim O(\varepsilon^N).
\end{align*}
 Repeating this process inductively yields all Sobolev estimates.
\end{enumerate}
 \end{proof}
\begin{remark}The proof of Proposition \ref{table_theorem} shows that any function $F\in\Sch(\RR)$ with $F(0)=1$ gives an embedding from $\maD'(M)$ to $\maG(M)$, but such an embedding will in general fail to  preserve  multiplication on smooth functions. An example is provided by the solution operator to the heat equation $e^{-\ep\Delta}$.
\end{remark}
\section{A larger class of embeddings}
 Here we generate new embeddings from old ones.
\begin{lem}
If $T_{\varepsilon}$ is a special embedding and  $N_{\varepsilon}$ a negligible net in $E(\Psi^{-\infty}(M))$ then $T_{\varepsilon}+N_{\varepsilon}$ is a special embedding.
\end{lem}
\begin{proof}
The negligibility of $N_{\varepsilon}$ in particular implies that for all $k$,  the Sobolev operator norm $\|N_{\varepsilon}\|_k\sim O(\varepsilon^N)$ where $N$ is any  integer. Thus in all Sobolev estimates one obtains no contributions from $N_{\varepsilon}$.
\end{proof}
\begin{prop}\label{converge}
Let $T_{n,\varepsilon}$ be a special embedding for all positive integers $n$, such that
\begin{enumerate}[a{)}]
\item\label{limit} For any distribution $u$,  $T_{n,\varepsilon}u\rightarrow u$ in $\maD'(M)$  as $\varepsilon\rightarrow 0$.
\item $T_{n,\varepsilon}$ is a moderate net in $E(\Psi^{-\infty}(M))$ and class $\ip{T_{n,\varepsilon}}$ converges in sharp topology in $\maG_{\Psi^{-\infty}(M)}$ to $T=T_{\varepsilon}$.
\end{enumerate}
Then any representative of $T$ defines a special embedding.
\end{prop}
\begin{remark}
The condition \ref{converge} \ref{limit}) above holds when $T_{\varepsilon}\rightarrow 
\mathrm{Id}$ in weak-$\ast$ operator topology in every $\mathscr{B}(H^k(M))$.
\end{remark}
\begin{proof}
Since $T_{\varepsilon}$ is a moderate net in $\Psi^{-\infty}(M)$, for any $u\in H^s(M)$ $\exists N$ such that
\[\|T_{\varepsilon}u\|_k\leq \|T_{\varepsilon}\|_{k+s}\|u\|_s\simeq O(\varepsilon^N).\]
 Let $u\in H^s(M)$. Due to the convergence of $T_n\rightarrow T$ in sharp topology there is an $n$ such  the operator norm $\|T_{\varepsilon}-T_{n,\varepsilon}\|_{H^s(M)\rightarrow C(M)}<C\varepsilon^2$ for small enough $\varepsilon$. Therefore
\[\lim_{\varepsilon\rightarrow 0}(T_{\varepsilon}-T_{n,\varepsilon})u=0\; \textrm{in}\;\maC(M).\]
Thus by embedding of continuous functions into $\maD'(M)$ we get
\begin{align*}
\lim_{\varepsilon\rightarrow 0}T_{\varepsilon}u&=\lim_{\varepsilon\rightarrow 0}(T_{\varepsilon}-T_{n,\varepsilon})u+ \lim_{\varepsilon\rightarrow 0}T_{n,\varepsilon}u\\
&=u.
\end{align*}
The invariance of multiplication of smooth functions follows analogously.
\end{proof}
Now we proceed to enlarge our class of Schwartz functions. Let $F_n$  be a sequence  that converges to $F$ in $\Sch(\RR)$ fast enough that is for every semi-norm $\rho$ there exists a  constant $C_{\rho}>0$ such that,
\begin{eqnarray}\label{fast_convergence}
\rho(F_n-F)<\frac{C_{\rho}}{n}\quad \textrm{for large}\, n.
\end{eqnarray}
 Given an integer $l$ we set
\begin{eqnarray}\label{speed}
k:=k(l,\varepsilon):=\inf\{j\in\ZZ|j>\frac{1}{\varepsilon^l}\}. 
\end{eqnarray}
Then we define $F^{[l]}_{\varepsilon}:=F_k(\varepsilon x)$.
\begin{lem}\label{sharp}
The nets $F^{[l]}$ define moderate elements in $E_{\Sch(\RR)}$. Furthermore $F^{[l]}_{\varepsilon}\rightarrow F_{\varepsilon}$ 
in the sharp topology on $\maG_{\Sch(\RR)}$.
\end{lem}
\begin{proof}
Since $F_n$ converge to $F$ in $\Sch(\RR)$  it is a bounded set of Schwartz functions and hence the moderateness of $F^{[l]}$ is immediate.

Let $N$ be a given integer and given $\alpha,\beta$ we observe that for all $l>N+|\alpha|+|\beta|$ one has,
\begin{align*}
\|x^{\beta}\partial^{\alpha}(F^{[l]}_{\varepsilon}-F_{\varepsilon})\|_{\infty}&= \varepsilon^{|\alpha|-|\beta|}\|x^{\beta}(F^{(\alpha)}_{K(l,\varepsilon)}(\varepsilon x)-F^{(\alpha)}(\varepsilon x))\|_{\infty}\\
&\leq \frac{C_{\alpha,\beta}\varepsilon^{|\alpha|-|\beta|}}{K(l,\varepsilon)}\qquad \textrm{where}\, C_{\alpha,\beta}\,\textrm{are given by \eqref{fast_convergence}}\\
&\leq O(\varepsilon^{l+|\alpha|-|\beta|})
\end{align*}
 Thus $F_{\varepsilon}^{[l]}\rightarrow F_{\varepsilon}$ in sharp topology.
\end{proof}
\begin{remark}
Given a net of  equi-continuous injection maps $i_{\varepsilon}:X\rightarrow X$ so that $x\rightarrow i_{\varepsilon}(x)$ defines an embedding of $X$ into its generalized space $\maG_X$ one could by the above recipe construct an sequence $Y^{[l]}_{\varepsilon}$   in $\maG_X$ from a rapidly converging sequence $x_n\rightarrow x$ such that $Y^{[l]}_{\varepsilon}\rightarrow i_{\varepsilon}(x)$ in sharp topology. 
\end{remark}
Now let $F_n\rightarrow F$ be a sequence in $\Sch(\RR)$ such that $F_n\equiv 1$ in an interval $(-t_n,t_n)$. We define analogous to \eqref{Weyl} for our nets $F^{[l]}_{\varepsilon}$;
\[N_l(\varepsilon):=\sup\{k|\lambda_k\leq \frac{t_{K(l,\varepsilon)}}{\varepsilon}\}.\]
\begin{lem}\label{Weyl2}
Let the sequence of intervals $t_n$ be such that for any $\alpha>0$ there exists $N_{\alpha}$ such that $t_n>n^{\!-\!\alpha}$ for all $n>N_{\alpha}.$ Then there exist $C>0,\tau>0$ with $N_l(\varepsilon)>C\varepsilon^{-\tau}$ for all $l$ and small enough $\varepsilon$.
\end{lem}
\begin{proof}
First we note that from (\ref{speed}) $K(l,\varepsilon)>\varepsilon^{-l}$. Therefore from the hypothesis on the $t_n$`s the corresponding intervals $t_{K(l,\varepsilon)}> \varepsilon^{l\alpha}$ for any $\alpha>0$. In particular $\frac{ t_{K(l,\varepsilon)}}{\varepsilon}\geq\varepsilon^{-1+\alpha}$. And the result follows again just as before by Weyl's asymptotic formula for $\lambda_k$`s. 
\end{proof}
\begin{definition}\label{admissible}
We say a Schwartz function $F\in \Sch(\RR)$ is admissible (or is almost $1$ near $0$) if there exists a sequence of Schwartz functions  $F_n$ with the following properties:
\begin{enumerate}[a{)}]
\item $F_n$ converges to  $F$ fast enough, that is for every continuous semi-norm $\rho$ on $\Sch(\RR)$  there exists a  constant $C_{\rho}>0$ such that,
\[\rho(F_n-F)<\frac{C_{\rho}}{n}\quad \textrm{for large}\, n.\]
\item There exists  $t_n$ such that $F_n\equiv 1$ on an interval $(-t_n,t_n)$ and such that for any $\alpha>0$ there exists $N_{\alpha}$ such that $t_n>n^{\!-\!\alpha}$ for all $n>N_{\alpha}.$ 
\end{enumerate}
\end{definition}

It is immediate that if a Schwartz function $F$ is admissible then $F(0)=1$ and $\partial^{\alpha} F(0)=0$
for all $\alpha$.

Before we state our main result we need the following lemma.
\begin{lem} Let $F_n\rightarrow F$ be a sequence in $\Sch(\RR)$ with $F_n\equiv 1$ 
in an interval $(-t_n,t_n)$ such that for any $\alpha>0$ there exists $N_{\alpha}$ such that $t_n< n^{\!-\!\alpha}$ for all $n>N_{\alpha}.$  Then the nets
\[F^{[l]}_{\varepsilon}(\Delta):=F_k(\varepsilon \Delta)\,\,\,\textrm{ where  }\, k:=K(l,\varepsilon):=\inf\{j\in\ZZ|j>\frac{1}{\varepsilon^l}\}\]
define special embeddings for all $l$.
\end{lem}
\begin{proof}
This is a consequence of  Lemma \ref{Weyl2}. Firstly since  each $F^{[l]}_{\varepsilon}(x)$ is a moderate net in $\Sch(\RR)$ so is $F^{[l]}_{\varepsilon}(\Delta)$. In particular $F^{[l]}_{\varepsilon}(\Delta)u$ is a moderate net in $\maG(X)$ for any distribution $u$.

Just like in  Proposition \ref{table_theorem} for  any $u\in\maD'(M)$ one checks that $\lim_{\varepsilon\rightarrow 0}F_{\varepsilon}^{[l]}(\Delta)u \rightarrow u$ therefore it maps $\maD'(M)\rightarrow \maG(M)$ injectively. First let us assume that  $u$ be in $L^2(M)$ so we may write a Fourier expansion $u=\sum a_n\phi_n$. 
Since by Lemma \ref{Weyl2} $N_l(\varepsilon)\rightarrow \infty$ for any integer $k$, there exist  $\varepsilon_0$ such that  $\varepsilon<\frac{t_{K(l,\varepsilon)}}	{\lambda_k}$ for all $0<\varepsilon<\varepsilon_0$. Now $u-F_{\varepsilon}^{[l]}(\Delta)(u)\leq\sum_{n>k}a_n(F_{\varepsilon}^{[l]}(\lambda_n)-1)\phi_n$ which tends to zero in $L^2(M)$ as $k\rightarrow \infty$ since all $F^{[l]}_{\varepsilon}(\lambda_n)$ are uniformally bounded. For $u$ in any other Sobolev space $H^s(M)$ we notice that the above argument can be applied to $(1+\Delta)^{-\frac{s}{2}}u$ and that these operators commute with $F^{[l]}_{\varepsilon}(\Delta)$.

The multiplicativity on smooth functions follows by a similar argument as before. We prove below the $L^2(M)$ estimate and all other estimates follow from it using  G\r{a}rding's inequality.

Let the Fourier expansion of a smooth function $f\in\smooth(M)$ be,
 \[f=\sum a_k\phi_k\]
then $(a_k)\in\Sch(\NN)$. Given any integer $N$, we can find a $k$ such that for all $l>k$
\[\sum_{n>l}|a_n|^2<\frac{1}{l^{\frac{N}{\tau}}},\]
where $\tau>0$ is chosen so that $N_l(\varepsilon)>(\varepsilon^{-\tau})$. Now choose $\varepsilon_0$ such that $\lambda_k<N_l(\varepsilon) $ for all $0<\varepsilon<\varepsilon_0$. Then for any $\varepsilon<\varepsilon_0$,
\[F^{[l]}_{\varepsilon}(\Delta)f-f=\sum_{n>N_l(\varepsilon)}a_n(F_{\varepsilon}(\lambda_n)-1)\phi_n.\]
 Since $F^{[l]}_{\varepsilon}(\lambda_n)$ are all uniformly bounded, 
 \begin{align*}
 \|F^{[l]}_{\varepsilon}(\Delta)f-f\|_{L^2(M)}^2&\leq C\sum_{n>N_l(\varepsilon)}|a_n|^2\\
 &\leq \frac{1}{N_l(\varepsilon)^{\frac{N}{\tau}}}\sim O(\varepsilon^N).
\end{align*}
\end{proof}
\begin{prop}\label{more_embeddings}
Let $F$ be an admissible Schwartz function. Then $F_{\varepsilon}(\Delta)$ is a special embedding.
\end{prop}
\begin{proof}
Since $F$ is admissible there is a sequence $F_n$ of Schwartz function approximating it in $\Sch(\RR)$ satisfying Definition \ref{admissible}.
Then by Lemma \ref{sharp} the sequence $F^{[l]}_{\varepsilon}$ defines a moderate net in $\maG_{\Sch(\RR)}$ converging in sharp topology to $F_{\varepsilon}$. Thus by applying $\Delta_*$, we get that $F^{[l]}_{\varepsilon}(\Delta)$ is a moderate net of smoothing operators that converges to $F_{\varepsilon}(\Delta)$. The result now follows in view of Proposition \ref{converge}  and the fact that each $F^{[l]}_{\varepsilon}(\Delta)$ is a special embedding. 
\end{proof}
\section{Invariance properties}
\subsection{Invariant operators}
We shall note how certain operators on the distributions behave  after they have been embedded into a generalized function algebra.
\begin{definition}
Let $D:\smooth(M)\rightarrow \smooth(M)$ be an operator that extends continuously to an operator $\tilde{D}:\maD'(M)\rightarrow \maD'(M)$. We say that $D$ is an invariant operator with respect to an embedding $T_{\varepsilon}$, if $\forall u\in \maD'(M)$,
\[(T_{\varepsilon}\tilde{D}u-DT_{\varepsilon}u)\in N_{\smooth(M)}.\]
\end{definition}
Let $\phi:M\rightarrow M$ be a diffeomorphism  then $\phi$ acts on $\smooth(M)$ by pull-back. The pull back extends on $\maD'(M)$ as:
\[\ip{\phi^*u,f}=\ip{u,\left.\phi^{\!-\!1}\right.^*f}.\]
\begin{prop}
If $\phi:M\rightarrow M$ is an isometry on a Riemannian manifold then $\phi$  is an invariant operator  with respect to the embedding $f_{\varepsilon}(\Delta)$, for all admissible Schwartz functions $f$ and the Laplace operator $\Delta$ associated with the metric.
\end{prop}
\begin{proof}
Since $\phi$ is an isometry $\phi_*(\Delta)=\Delta$. Therefore $\phi$ preserves the eigenspaces for all operators of the form $F(\Delta)$ and hence $\phi_*(F(\Delta)=F(\Delta)$. Therefore the result follows.
\end{proof}
 Of course all pseudodifferential operators which commute with the Laplace $\Delta$  also are preserved under embeddings obtained from $\Delta$.

\subsection{Support and singular support}
We recall that $\maG^{\infty}(M)$ is a sub-algebra of $\maG(M)$  so that $u\in\maG^{\infty}(M)$ iff there exists an $N$ such that for any continuous semi-norm $\rho$ on $\smooth(M)$,  $\rho(u)\sim O(\varepsilon^N)$. That is, the moderateness estimates in $\maG^{\infty}(M)$ are independent of the semi-norms. As introduced in \cite{MObook}
 on $\RR^n$ the convolution embedding has the property that  $\iota_\rho({\mathcal D}')\cap \maG^{\infty} = {\mathcal C}^\infty$. We shall need a generalization of this result 
 this result to our geometrical embeddings:

\begin{prop}\label{regularity}
Let $F$ be a Schwartz function such that $F\equiv 1$ on $(-t,t)$. Then the embedding of distributions  under $T_{\varepsilon}=F_{\varepsilon}(\Delta)$ satisfies 
 \[T_{\varepsilon}\maD'(M)\cap \maG^{\infty}(M)=\smooth(M)\,.\]
%\[WF(u)=WF_g(F_{\varepsilon}(\Delta)u).\]
\end{prop}
\begin{proof}
 It is clear that any smooth function embeds into $\maG^{\infty}(M)$. 
We shall  prove that any non-smooth distribution is not in $\maG^{\infty}(M)$.
  To start with we make the following claim:
\begin{lem}
Let $m= \operatorname{dim}(M)$. Let $u$ be a non-smooth distribution then there exists an $s$ such that  $u$ is in $H^s(M)$ and  $u\not\in H^t(M)$ for any $t>s$. 
 Then
%\begin{enumerate}
%\item For any $\delta>0$  there exists a $C$ such that as $\ep\rightarrow 0$.
%\[\|F_{\varepsilon}(\Delta)u\|_{L^2(M)}\leq C\ep^{-\frac{1}{2}-\frac{s}{m}-\delta}.\]
%Or the net $\|F_{\varepsilon}(\Delta)u\|_{L^2(M)}$ is $O(\ep^{-\frac{1}{2}-\frac{s}{m}-\delta})$.
%\item
 for any $\delta>0$ and for any $C>0$ there exists a sequence $\ep_n\rightarrow 0$ such that
\[\|F_{\varepsilon_n}(\Delta)u\|_{L^2(M)}> C\ep_n^{\frac{m}{2}+s+\delta}\,,
\]
i.e., the net $\|F_{\varepsilon}(\Delta)u\|_{L^2(M)}$ is not $O(\ep^{\frac{m}{2}+s+\delta}).$
%\end{enumerate}
\end{lem}
\begin{proof}
 Let $\phi_k$ be the $L^2(M)$ spectral basis for Laplace with nondecreasing eigenvalues. Since $(1+\Delta)^{\frac{s}{2}}u\in\ L^2(M)$
   there exists $b_k\in \ell^2$ such that,
\begin{align*}
(1+\Delta)^{\frac{s}{2}}u&=\sum b_k\phi_k\\
F_{\varepsilon}(\Delta)u&=\sum F_{\varepsilon}(\lambda_k)(1+\lambda_k)^{-\frac{s}{2}}b_k\phi_k.
\end{align*}
  Let $a_k=(1+\lambda)^{-\frac{s}{2}}b_k$. By Weyl's theorem $\lambda_k\sim Ck^{\frac{2}{m}}$. Therefore with  our assumptions on $u$ and on $s$ we have:
\begin{align}
a_k(1+\lambda_k)^{\frac{s}{2}}\in\ell^2&\Rightarrow a_kk^{\frac{s}{m}}\in \ell^2\notag\\
\textrm{for}\,\,t>s\,\,a_k(1+\lambda_k)^{\frac{t}{2}}\not\in \ell^2 &\Rightarrow a_kk^{\frac{t}{m}}\not\in \ell^2\label{unbounded}
\end{align}
 Hence by  using (\ref{unbounded}), given a constant $ C>0$ we get a sequence $k_n$ such that $|a_{k_n}|>Ck^{-\frac{s}{m}-\frac{1}{2}-\delta}$. (This is seen by comparing the divergent series $|a_k|^2k^{\frac{2s}{m}+\delta}$ with the convergent series $Ck^{-1-\delta}$.)

 Now pick $\ep_n$ such that 
\begin{enumerate}[(a)]
\item if $-1-\frac{2s}{m}-\delta>0$ then

\[
N_{\ep_n+\tau}<k_n\leq N_{\ep}\,\textrm{ for all}\,\tau>0.\]
\item
Else if $ -1-\frac{2s}{m}-\delta<0$
\[N_{\ep_n}<k_n\leq N_{\ep-\tau}\,\textrm{ for all}\,\tau>0
\]
\end{enumerate}
Since $k_n\rightarrow \infty$ therefore $\ep_n\rightarrow 0$ and
 
\begin{align*}
\|F_{\varepsilon_n}(\Delta)u\|_{L^2(M)}^2&>\sum_{n\leq N_{\varepsilon_n}} |a_nF_{\varepsilon_n}(\lambda_n)|^2\\
&=\sum_{n\leq N_{\varepsilon_n}} |a_n|^2\\
&>CN(\varepsilon_n)^{-1-\frac{2s}{m}-\delta}>C\varepsilon_n^{\frac{m}{2}+s+\delta}
\end{align*}
by Lemma \ref{asymptotic}. 
\end{proof}
Thus the $L^2(M)$-estimate of $u$ has an asymptotic that depends on $s$.
Any $k$--th Sobolev estimate for $u$  is in fact $L^2(M)$ estimate for 
$(1+\Delta)^{\frac{k}{2}}u\in H^{s-k}(M)$. From the lemma above and the fact that $u\not\in H^t(M)$ 
implies  $(1+\Delta)^{\frac{k}{2}}u\not\in H^{t-k}(M)$ we obtain $F_{\varepsilon}(\Delta)u\not\in\maG^{\infty}(M)$. 
\end{proof}
 Since $\maG(\,\_\,)$ is in fact a sheaf of algebras (cf.\ \cite{RD,book}), there is a well-defined 
notion of support of its elements.
\begin{lem}\label{support}
The embeddings $F_{\varepsilon}(\Delta)$ preserve the support of distributions, that is 
$\operatorname{supp}(u)=\operatorname{supp}(F_{\varepsilon}(\Delta)u)$ for any distribution $u$  
with compact support.
\end{lem}
\begin{proof}
This is a consequence of finite propagation speed for operators with kernels supported near the diagonal.
Given $\delta>0$ we pick a cutoff function 
$\phi_{\delta}(x)=\begin{cases} 1&|x]<\frac{\delta}{2}\\0&|x|>\delta
\end{cases}$\\
 Then $(1-\phi_{\delta}(x))F_{\varepsilon}(x)$  is a negligible net in $N_{\Sch(\RR)}$  and therefor  the classes $F_{\varepsilon}(\Delta)u$ is the same as $\phi_{\delta}(\Delta))F_{\varepsilon}(\Delta)u$. Since the operators $\phi_{\delta}(\Delta))F_{\varepsilon}(\Delta)$ are supported in a $\delta$ neighborhood of the diagonal the support of $\phi_{\delta}(x))F_{\varepsilon}(\Delta)u$  is in $2\delta$ of $\operatorname{supp}u$. Since $\delta$ was arbitrary  we have the desired result.
\end{proof}
The singular support for generalized functions can be defined in terms of $\maG^{\infty}(M)$ as
\begin{multline*}
x\not\ni\operatorname{singsupp}u\Leftrightarrow \exists\,\phi\in\smooth(M)\, \textrm{ s.t.\ }
\phi(x)\neq0\,\textrm{and}\,\phi(x)u\in\maG^{\infty}(M)
\end{multline*}
 This is  analogous to the definition of singular support for a distribution. Writing a distribution $u$ as 
\[u=\phi u+(1-\phi)u,\]
 we note that  by Lemma \ref{regularity}
\[F_{\ep}(\Delta)\phi u\in\maG^{\infty}(M)\Longleftrightarrow \phi u\in\smooth(M).\]
Moreover, by Lemma \ref{support} if $\phi$ is supported near $x$ then $F_{\ep}(\Delta)(1-\phi) u$ is supported away from $x$. Therefore, 
\begin{cor}
The embeddings $F_{\varepsilon}(\Delta)$ preserve the singular support of distributions, that is 
$\operatorname{singsupp}(u)=\operatorname{singsupp}(F_{\varepsilon}(\Delta)u)$ for any distribution $u$.
\end{cor}
Let $P$ be a (classical) pseudodifferential operator of order $0$ with $\sigma(P)$  its principal symbol. The characteristic set of $P$  is $char(P)=\sigma(P)^{-1}(0)\subseteq T^*M\setminus 0$. 
We follow the results in \cite{garetto_hoermann} to define the generalized wavefront set of $u\in\maG(M)$ to be
\[WF_g(u):=\bigcap_{Pu\in\maG^{\infty}(M)} \operatorname{char}(P).\]
 The preservation of wavefront set under our embeddings shall be considered elsewhere.

\section{Sections of vector bundles}
 We begin by reviewing the construction of generalized sections of a vector bundle (\cite{ndg,gprg,book}) 
 and then construct embeddings of distributional sections into generalized sections preserving various properties of interest.
Let $E\rightarrow M$ be a (complex) vector bundle over $M$. We fix a fiberwise inner product on $E$. Let $\Gamma^\infty(E)=\Gamma^\infty(M,\,E)$ be the space of smooth sections of $E$ and $L^2(E)$ its closure with respect to the inner product
\[\ip{s_1,s_2}:=\int_M\ip{s_1(x),\overline{s_2(x)}}dx,\qquad s_1,s_2\in\Gamma^\infty(E).\]
We shall denote the norm on this Hilbert space by $\|\,\|_{L^2(E)}$.

Let  $D:\gp{E}\rightarrow \gp{E}$  be an order one positive elliptic operator  on $E$. 
Then the Sobolev spaces $H^t(E)$ as the completion of $\gp{E}$ with respect to the norm:
\[\|s\|_t=\|D^t s\|_{L^2(E)}.\]

$\Gamma^\infty(E)$ is a locally convex space with respect to the family of norms $\|\,\|_t$. The 
$\tilde{\CC}$-module associated with the Fr\'echet space $\gp{E}$ will be denoted by  $\maG_E(M)$ or $\maG_E$ when no confusion can arise. We note the following obvious properties.
\begin{lem}
Let $E$ and $F$ be vector bundles over $M$ then
\begin{enumerate}
\item \label{mod}Each $\maG_E(M)$ is a module over $\maG(M)$.
\item \label{whit}As $\maG(M)$ modules the generalized sections of the Whitney sum are given by 
$\maG_{E\oplus F}(M)=\maG_E(M)\oplus \maG_F(M)$ .
\item\label{tensor} There is an inclusion of $\maG(M)$ modules \[i:\maG_E(M)\otimes_{\maG(M)} \maG_F(M)\rightarrow \maG_{E\otimes F}(M).\]
\item \label{opt} A pseudodifferential operator $P:\gp{E}\rightarrow \gp{F}$ induces a  map $\hat{P}: \maG_E(M)\rightarrow \maG_F(M)$.
\end{enumerate}
\end{lem}
\begin{proof}
Part (\ref{mod}) and (\ref{whit}) are obvious. For (\ref{tensor})  we note that $\ip{x_{\varepsilon}}\otimes \ip{y_{\varepsilon}}\rightarrow \ip{x_{\varepsilon}\otimes y_{\varepsilon}}$ is the required injection, as $\gp{E}\otimes_{\smooth(M)} \gp{F}\simeq \gp{E\otimes F}$ are isomorphic.(\ref{opt}) is a special case of Lemma \ref{lcs}.
\end{proof}
 In particular the de Rham differential $d:\Omega^*(M)\rightarrow \Omega^*(M)$ induces a map $\hat{d}$ on the generalized differential forms (cf.\ \cite{ndg,connections}).

Let $D$ be a positive elliptic operator on sections of a vector bundle $E$ such that the order of $D$ is $\geq 1$. Then Weyl's estimate  holds for the spectrum of $D$, namely there is a spectral decomposition of $L^2(E)$ in terms of eigenvectors of $E$  and the corresponding eigenvalues satisfy asymptotic estimates of rational order. This in particular implies that we can use the Schwartz functional calculus. For any $F\in\Sch(\RR)$ the operator $F(D)$ is a smoothing operator on $E$. Such an operator has kernel in the big homomorphism bundle
\[\ker(f(D))\in \gp{M\times M:\pi_L^*E\otimes \pi_R^*E'\otimes \pi_R^*\Omega_M},\]
 where $\Omega_M$ is the 1-density bundle on $M$ and $E'$ is the dual bundle. We  could again  fix a Riemannian density for once and  normalize the kernel. The following can be proven in exactly the same way as Proposition \ref{more_embeddings}.
\begin{prop}
If $F\in\Sch(\RR)$ is an admissible Schwartz function then the  net  of operators $F_{\varepsilon}(D):\maD'(E)\rightarrow \maG_E(M)$  is an embedding that coincides  with the constant  embedding of smooth sections $\gp{E}$. In particular the following diagram commutes with respect to the two  module actions,
\begin{center}
$ \xymatrix{ {\maG(M)\times \maG_E(M)}\ar[r]&
{\maG_E(M)}\\ 
{\smooth(M)\times \gp{E}}\ar@<3ex>[u]^{F_{\varepsilon}(\Delta)}\ar@<-3ex>[u]_{F_{\varepsilon}(D)}\ar[r]&\gp{E}\ar[u]^{F_{\varepsilon}(D)}}$.
\end{center}
\end{prop}
 We can now for example define a generalized connection on a vector bundle $E$ as a $\tilde{\CC}$-linear map
\[\nabla:\maG_{E\,\otimes\, \Omega^*(M)}\rightarrow \maG_{E\,\otimes\, \Omega^*(M)}\]
 such that
\[\nabla(s\otimes\alpha)=\nabla(s)\otimes \alpha+(-1)^ks\otimes \hat{d}(\alpha),\qquad s\in\maG_{E},\, \alpha\in\maG_{\Omega^k(M)}.\] 
Of course any usual connection on $E$ gives rise to a generalized connections. Generalized connections 
are studied in \cite{connections} 
%\bibliographystyle{siam}
%\bibliography{diana}

\end{document}